\documentclass[12pt,reqno]{amsart}

\usepackage{amsmath,amsthm,amsfonts,amscd,amssymb,euscript,enumerate}
\numberwithin{equation}{section} 
\newcommand{\dbar}{\ensuremath{\bar \partial}}

\newcommand{\C}{\ensuremath{{\mathbb C}}}
\newcommand{\pj}{\ensuremath{{\mathbb P}}}

\newcommand{\N}{\ensuremath{{\mathbb N}}}

\newcommand{\smooth}{\ensuremath{C^{\infty}}}

\newcommand{\I}{\ensuremath{\mathcal I}}

\newcommand{\U}{\ensuremath{\mathcal U}}

\newcommand{\oka}{\ensuremath{\mathcal O}}

\newcommand{\cv}{\ensuremath{\mathcal C}}

\newtheorem{Lemma}{Lemma}
\newtheorem{Theorem}[Lemma]{Theorem}
\newtheorem{MainTheorem}[Lemma]{Main Theorem}
\newtheorem{Proposition}[Lemma]{Proposition}
\newtheorem{Corollary}[Lemma]{Corollary}
\newtheorem*{Theorem*}{Theorem}
\newtheorem*{thm1.1}{Theorem 1.1}
\newtheorem*{thm1.2}{Theorem 1.2}
\newtheorem*{thm2.9'}{Theorem 2.9'}
\newtheorem*{prop4.1}{Proposition 4.1}
\newtheorem*{prop4.2}{Proposition 4.2}
\newtheorem*{prop4.4}{Proposition 4.4}
\newtheorem*{prop4.5}{Proposition 4.5}
\theoremstyle{remark}
\newtheorem{Remark}[Lemma]{Remark}

\newtheorem*{Remark*}{Remark}
\theoremstyle{definition}
\newtheorem{Definition}[Lemma]{Definition}
\newtheorem*{Definition*}{Definition}

\numberwithin{Lemma}{section}

\numberwithin{equation}{section}

\makeatletter
\@namedef{subjclassname@2020}{%
  \textup{2020} Mathematics Subject Classification}
\makeatother

\begin{document}
\title[Relating Catlin and D'Angelo $q$-types]{Relating Catlin and D'Angelo $q$-types}
\author{Vasile Brinzanescu}

\address{Simion Stoilow Institute of Mathematics of the Romanian Academy, Research unit 3, 21 Calea Grivitei Street, 010702 Bucharest, Romania}

\email{Vasile.Brinzanescu@imar.ro}

\author{Andreea C. Nicoara}

\address{School of Mathematics, Trinity College Dublin, Dublin 2, Ireland}

\email{anicoara@maths.tcd.ie}



\begin{abstract}
We clarify the relationship between the two most standard measurements of the order of contact of q-dimensional complex varieties with a real hypersurface, the Catlin and D'Angelo $q$-types, by showing that the former equals the generic value of the normalized order of contact measured along curves whose infimum is by definition the D'Angelo $q$-type.  
\end{abstract}

\maketitle

\tableofcontents

\section{Introduction}
\label{intro}

The purpose of this note is finishing the work initiated in \cite{bazilandreea} as far as elucidating the relationship between Catlin and D'Angelo $q$-types. D'Angelo's finite $q$-type defined in 1982 in \cite{opendangelo} is by far the most standard finite type notion in several complex variables. The Kohn Conjecture, one of the most famous open problems in complex analysis, if proven would say that finite D'Angelo $q$-type at a certain point on the boundary of a smooth pseudoconvex domain in $\C^n$ ensures the termination of the Kohn algorithm defined in \cite{kohnacta} and thus the subellipticity of the $\dbar$-Neumann problem for $(p,q)$ forms in the neighborhood of that point. As of now, the only result on the subellipticity of the $\dbar$-Neumann problem for $(p,q)$ forms on smooth pseudoconvex domains in $\C^n$ is due to Catlin in  \cite{catlinsubell} and provides a lower bound on the subelliptic gain for the $\dbar$-Neumann problem in terms of his own notion of $q$-type. As a result, relating the D'Angelo and Catlin $q$-types effectively is paramount so that if and when the Kohn Conjecture is proven with a lower bound for subellipticity in place, Catlin's bound and the bound obtained via the Kohn algorithm can be compared.

As the results in \cite{bazilandreea} are proven with respect to the generic value, we introduce the following definition; see also \cite{bazilandreeaerr}:

\begin{Definition}
Let $2 \leq q \leq n.$ If $\I$ is an ideal in $\oka_{x_0},$
$$\tilde \Delta_q(\I, x_0)={\text gen.val}_{\{w_1, \dots, w_{q-1} \}} \:\:\ \Delta_1\Big( (\I, w_1, \dots, w_{q-1}), x_0 \Big),$$ where the generic value is taken over all non-degenerate sets $\{w_1, \dots, w_{q-1} \}$ of linear forms in $\oka_{x_0},$ $(\I, w_1, \dots, w_{q-1})$ is the ideal in $\oka_{x_0}$ generated by $\I, w_1, \dots, w_{q-1},$ and $\Delta_1$ is the D'Angelo $1$-type. Likewise, if $M$ is a real hypersurface in $\C^n$ and $x_0 \in M,$ $$\tilde \Delta_q (M, x_0) ={\text gen.val}_{\{w_1, \dots, w_{q-1} \}} \:\:\ \Delta_1\Big( (\I(M), w_1, \dots, w_{q-1}), x_0 \Big),$$ where $(\I(M), w_1, \dots, w_{q-1})$ is the ideal in $C^\infty_{x_0}$ generated by all smooth functions $\I(M)$ vanishing on $M$ along with $w_1, \dots, w_{q-1}.$ 
\end{Definition}

Let $D_q$ denote the Catlin $q$-type, and let $\Delta_q$ be the D'Angelo $q$-type. The two main theorems in \cite{bazilandreea} can now be stated with respect to $\tilde \Delta_q$ as follows:

\begin{Theorem}[Theorem 1.1 \cite{bazilandreea},\cite{bazilandreeaerr}]
\label{bazilandreeathm1.1}
Let $\I$ be an ideal of germs of holomorphic functions at $x_0,$ then for $1 \leq q \leq n$ $$D_q(\I, x_0) \leq \tilde \Delta_q(\I, x_0) \leq \left(D_q(\I, x_0)\right)^{n-q+1}.$$
\end{Theorem}

\begin{Theorem}[Theorem 1.2 \cite{bazilandreea},\cite{bazilandreeaerr}]
\label{bazilandreeathm1.2}
Let $\Omega$ in $\C^n$ be a domain with $\smooth$
boundary. Let $x_0 \in b \Omega$ be a point on the boundary of
the domain, and let $1 \leq q <n.$
\begin{enumerate}[(i)]
\item $D_q(b \Omega, x_0) \leq \tilde \Delta_q(b \Omega, x_0);$
\item If $\tilde \Delta_q(b \Omega, x_0)<\infty$ and the domain is $q$-positive at $x_0$ (the $q$ version of D'Angelo's property P), then $$\tilde \Delta_q(b \Omega, x_0) \leq 2 \left(\frac{ D_q(b \Omega, x_0)}{2} \right)^{n-q}.$$
\end{enumerate}
 In particular, if $b \Omega$ is pseudoconvex at $x_0$ and $\tilde \Delta_q(b \Omega, x_0)<\infty,$ then $$D_q(b \Omega, x_0) \leq \tilde \Delta_q(b \Omega, x_0) \leq 2 \left(\frac{ D_q(b \Omega, x_0)}{2} \right)^{n-q}.$$

\end{Theorem}

It is an easy consequence of D'Angelo's work in \cite{opendangelo} that $\Delta_q$ and $\tilde \Delta_q$ are simultaneously finite for ideals of germs of holomorphic functions:

\medskip
\begin{Proposition}
\label{normalandtilde}
Let $\I$ be any ideal in $\oka_{x_0}.$ For any $2 \leq q \leq n,$ $$\Delta_q(\I, x_0) \leq \tilde \Delta_q(\I, x_0) \leq \left( \Delta_q(\I, x_0) \right)^{n-q+1}.$$
\end{Proposition}

The counterpart of this result for points on the boundary of a domain is slightly more difficult and requires the assumption of $q$-positivity:

\medskip
\begin{Proposition}
\label{normalandtildebdry}
Let $\Omega$ in $\C^n$ be a domain with $\smooth$
boundary. Let $x_0 \in b \Omega$ be a point on the boundary of
the domain, and let $2 \leq q <n.$ If the domain is $q$-positive at $x_0,$ then
$$\Delta_q(b \Omega, x_0) \leq \tilde \Delta_q(b \Omega, x_0) \leq 2 \left( \Delta_q(b \Omega, x_0) \right)^{n-q}.$$ If the domain is pseudoconvex at $x_0,$ then it is $q$-positive, so the same inequality holds.
\end{Proposition}

We shall prove Propositions~\ref{normalandtilde} and ~\ref{normalandtildebdry} in Section~\ref{pfmainthm}. We note here at for $q=1$ the notions of $\Delta_q,$ $\tilde \Delta_q,$ and $D_q$ coincide, a fact obvious from their definitions.

\smallskip Martino Fassina produced an example in \cite{fassina} when the Catlin and D'Angelo $q$-types are not equal to each other, answering a question that has been open since 1987 when \cite{catlinsubell} was published. In that example, $$\Delta_2(\I,0)= 3 < 4 = \tilde \Delta_2(\I,0) = D_2(\I,0).$$ Furthermore, in the same paper \cite{fassina}, Fassina proved that given any positive integer, an ideal of holomorphic germs can be constructed so that the difference between the Catlin $q$-type and the D'Angelo $q$-type of that ideal is larger than the given integer. Fassina's work thus highlights the importance of finding an effective relationship between $\Delta_q$ and $D_q.$

The main result of this paper is the following characterization of the Catlin $q$-type $D_q$ that clarifies completely how it relates to the D'Angelo $q$-type:

\begin{MainTheorem}
\label{mainthm}

\begin{enumerate}[(i)]
\item Let $\I$ be an ideal of germs of holomorphic functions at $x_0,$ and let $2\leq q \leq n,$ then $$D_q(\I, x_0) = \tilde \Delta_q(\I, x_0).$$
\item Let $\Omega$ in $\C^n$ be a domain with $\smooth$ boundary. Let $x_0 \in b \Omega$ be a point on the boundary of the domain, and let $2 \leq q <n.$ Then $$D_q(b \Omega, x_0) = \tilde \Delta_q(b \Omega, x_0).$$
\end{enumerate}
\end{MainTheorem}

\smallskip
\begin{Corollary}
\label{maincor}
\begin{enumerate}[(i)]
\item Let $\I$ be an ideal of germs of holomorphic functions at $x_0,$ and let $2\leq q \leq n,$ then $$\Delta_q(\I, x_0) \leq D_q(\I, x_0) \leq \left( \Delta_q(\I, x_0) \right)^{n-q+1}.$$
\item Let $\Omega$ in $\C^n$ be a domain with $\smooth$ boundary. Let $x_0 \in b \Omega$ be a point on the boundary of the domain, let $2 \leq q <n,$ and assume that the domain is $q$-positive at $x_0.$ Then $$\Delta_q(b \Omega, x_0) \leq D_q(b \Omega, x_0) \leq 2 \left( \Delta_q(b \Omega, x_0) \right)^{n-q}.$$ In particular, if the domain is pseudoconvex at $x_0,$ then the same inequality holds.
\end{enumerate}
\end{Corollary}

In \cite{catlinsubell} Catlin defined his $q$-type $D_q$ by starting with the germ of a $q$-dimensional variety $V^q$ and constructing an open set in the Grassmannian $G^{n-q+1}$ of all $(n-q+1)$-dimensional complex linear subspaces through $x_0$ in $\C^n$ that is specific to a given $V^q.$ His construction is so delicate because the aim is to obtain the same number of curves in the intersection of $V^q$ with each $(n-q+1)$-dimensional complex linear subspace in this open set and the same maximal normalized order of contact measured along the intersection curves. As a result, proving the equality of $D_q$ with $\tilde \Delta_q$ is trickier than it seems. Given any curve and any subspace $W,$ we have to show the existence of a cylinder variety with the subspace $W$ as the directrix subspace along the curve at $x_0.$ This construction allows us to construct a $q$-dimensional variety $V^q$ starting with any curve whose open set in the Grassmanian $G^{n-q+1}$ in Catlin's construction is well behaved.

The paper is organized as follows: We devote Section~\ref{definitionssection} to recalling the definitions of the Catlin and D'Angelo $q$-types and relevant results. We then prove Propositions~\ref{normalandtilde} and \ref{normalandtildebdry} as well as the Main Theorem, Theorem~\ref{mainthm}, and Corollary~\ref{maincor} in Section~\ref{pfmainthm}.

\smallskip
\noindent {\bf Acknowledgements}
The work of Vasile Brinzanescu was partially supported by a grant of the Ministry of Research and Innovation, CNCS - UEFISCDI, project number PN-III-P4-ID-PCE-2016-0030, within PNCDI III, and CNCS-UEFISCDI project PN-III-P4-ID-PCE-2020-0029, Syzygies, invariants and classification problems in algebraic geometry and topology.

\section{Catlin and D'Angelo $q$-types}
\label{definitionssection}

For the convenience of the reader, we first recall the definitions of the D'Angelo and Catlin $q$-types and the properties needed for proving the results in Section~\ref{intro}.

\medskip
\begin{Definition}
\label{idealtype}
Let $C^\infty_{x_0}$ be the ring of smooth germs at $x_0 \in \C^n,$ and let $\I$ be an ideal in $C^\infty_{x_0}$ or $\oka_{x_0}.$ $$\Delta_1(\I,x_0)   = \sup_{\varphi \in \cv(n,x_0)} \:\: \inf_{g \in \I} \:\:\frac{ {\text ord}_0 \, \varphi^* g}{{\text ord}_0 \, \varphi},$$ where $\cv(n,x_0)$ is the set of all germs of holomorphic curves $$\varphi: (U,0) \rightarrow (\C^n, x_0)$$ such that $\varphi(0)=x_0$ for $U$ is some neighborhood of the origin in $\C^1,$ ${\text ord}_0$ is the vanishing order at the origin, and  ${\text ord}_0 \, \varphi = \min_{1 \leq j \leq m} \, {\text ord}_0 \, \varphi_j.$
\end{Definition}

\medskip
\begin{Definition}
Let $2 \leq q \leq n.$ If $\I$ is an ideal in $\oka_{x_0},$ the D'Angelo $q$-type is given by
$$\Delta_q(\I, x_0)=\inf_{\{w_1, \dots, w_{q-1} \}} \:\:\ \Delta_1\Big( (\I, w_1, \dots, w_{q-1}), x_0 \Big),$$ where the infimum is taken over all non-degenerate sets $\{w_1, \dots, w_{q-1} \}$ of linear forms in $\oka_{x_0}$ and $(\I, w_1, \dots, w_{q-1})$ is the ideal in $\oka_{x_0}$ generated by $\I, w_1, \dots, w_{q-1}.$ Likewise, if $M$ is a real hypersurface in $\C^n,$ $x_0 \in M,$ and $2 \leq q < n,$ the D'Angelo $q$-type of the hypersurface $M$ is given by $$\Delta_q (M, x_0) =\inf_{\{w_1, \dots, w_{q-1} \}} \:\:\ \Delta_1\Big( (\I(M), w_1, \dots, w_{q-1}), x_0 \Big),$$ where $(\I(M), w_1, \dots, w_{q-1})$ is the ideal in $C^\infty_{x_0}$ generated by all smooth functions $\I(M)$ vanishing on $M$ along with $w_1, \dots, w_{q-1}.$ 
\end{Definition}

\medskip
\begin{Definition}
If $\I$ is an ideal in $\oka_{x_0},$ $$D(\I, x_0) = \dim_\C (\oka_{x_0} / \I).$$
\end{Definition}

\medskip
\begin{Proposition}[Proposition 2.8 \cite{bazilandreea}]
\label{multgeneric}
Let $\I$ be a proper ideal in $\oka_{x_0},$ and let $x_0 \in \C^n.$
\begin{equation*}
\inf_{\{w_1, \dots, w_{q-1} \}} D\Big( (\I, w_1, \dots, w_{q-1}), x_0 \Big)={\text gen.val}_{\{w_1, \dots, w_{q-1} \}} D\Big( (\I, w_1, \dots, w_{q-1}), x_0 \Big),
\end{equation*}
where $\{w_1, \dots, w_{q-1} \}$ is a non-degenerate set of linear forms in $\oka_{x_0},$ $(\I, w_1, \dots, w_{q-1})$ is the ideal in $\oka_{x_0}$ generated by $\I, w_1, \dots, w_{q-1},$ and the infimum and the generic value are both taken over all such non-degenerate sets $\{w_1, \dots, w_{q-1} \}$ of linear forms in $\oka_{x_0}.$ In other words, the infimum is achieved and equals the generic value.
\end{Proposition}

\medskip
\begin{Theorem}[D'Angelo, Theorem 2.7 \cite{opendangelo}]
\label{opendangelomain}
If $\I$ is an ideal in $\oka_{x_0}$ containing $q$ linearly independent linear forms $w_1, \dots, w_q,$ then 
\begin{equation}
\label{dangeloineq}
\Delta_1(\I, x_0) \leq D(\I, x_0) \leq \left( \Delta_1(\I,x_0) \right)^{n-q}.
\end{equation}
\end{Theorem}

\bigskip\noindent If $\Delta_q (M, x_0) = t < \infty,$ let $k = \lceil t \rceil$ be the ceiling of $t.$ By Proposition 14 from p.88 of \cite{dangeloftc}, 
\begin{equation}
\label{truncation}
\Delta_q (M, x_0) = \Delta_q (M_k, x_0),
\end{equation} 
where $M_k$ is real hypersurface defined by $r_k,$ the polynomial with the same $k$-jet at $x_0$ as the defining function $r$ of $M.$ By D'Angelo's polarization technique from \cite{opendangelo}, $$r_k = Re\{h\} + ||f||^2-||g||^2,$$ where $||f||^2= \sum_{j=1}^N |f_j|^2,$ $||g||^2= \sum_{j=1}^N |g_j|^2,$ and the functions $h, f_1, \dots, f_N, g_1, \dots, g_N$ are all holomorphic polynomials in $n$ variables. Let $\U(N)$ be the group of $N \times N$ unitary matrices. $\forall \, U \in \, \U(N)$ consider the ideal of holomorphic polynomials $\I(U, x_0)=(h, f-Ug)$ generated by $h$ and the $N$ components of $f-Ug,$ where $f=(f_1, \dots, f_N)$ and $g= (g_1, \dots, g_N).$ D'Angelo proved the following result:

\medskip
\begin{Theorem}[Theorem 14 \cite{dangeloftc}]
\label{dangeloqboundthm}
\begin{equation*}
\Delta_q(M, x_0) \leq 2   \sup_{U \in \, \U(N)} \:\: \inf_{\{w_1, \dots, w_{q-1} \}} \:\: D\Big( (\I(U,x_0), w_1, \dots, w_{q-1}), x_0 \Big) \leq 2 \Big(\Delta_q(M, x_0) \Big)^{n-q}
\end{equation*}
\end{Theorem}

\medskip\noindent Let us also note that $\U(N)$ is compact, so $\displaystyle \sup_{U \in \, \U(N)} \: D\Big( (\I(U,x_0), w_1, \dots, w_{q-1})$ is achieved due to the fact that the multiplicity of an ideal is upper semi-continuous:

\medskip
\begin{Proposition}[Part of Proposition II.5.3 \cite{tougeron}]
\label{duscprop}
Let $\I(\lambda)$ be an ideal in $\oka_{x_0}$ that depends continuously on $\lambda.$ Then $D\big(\I(\lambda), x_0\big)$ is an upper semi-continuous function of $\lambda.$ 
\end{Proposition}

\bigskip\noindent We recall from \cite{bazilandreea} the $q$ version of D'Angelo's property P, $q$-positivity, the hypothesis that appears in Theorem~\ref{bazilandreeathm1.2}:

\medskip
\begin{Definition}[Definition 2.14 \cite{bazilandreea}]
\label{qpropertyP}
Let $M$ be a real hypersurface of $\C^n,$ and let $x_0 \in M$ be such that $\Delta_q(M, x_0) <k.$ Let $j_{k,x_0} r = r_k = Re\{h\} + ||f||^2-||g||^2$ be a holomorphic decomposition at $x_0$ of the $k$-jet of the defining function $r$ of $M.$ We say that $M$ is $q$-positive at $x_0$ if for every holomorphic curve $\varphi \in \cv(n,x_0)$ for which $\varphi^* h$ vanishes and such that the image of $\varphi$ locally lies in the zero locus of a non-degenerate set of linear forms $\{w_1, \dots, w_{q-1} \}$ at $x_0,$ the following two conditions are satisfied:
\begin{enumerate}
\item[(i)] ${\text ord}_0 \, \varphi^* r$ is even, i.e. ${\text ord}_0 \, \varphi^* r = 2a,$ for some $a \in \N;$
\item[(ii)] $\displaystyle \left(\frac{d}{dt}\right)^a \left(\frac{d}{d \bar t}\right)^a \varphi^* r (0) \neq 0.$
\end{enumerate}
\end{Definition}

\medskip When $q=1$ this definition is exactly D'Angelo's property P with respect to which he proved the following result:

\medskip
\begin{Theorem}[Theorem 5.3 \cite{opendangelo}]
\label{1propertyP}
If $M$ satisfies property P at $x_0,$ then $$\Delta_1(M, x_0) = 2 \,  \sup_{U \in \, \U(N)} \:\: \Delta_1(\I(U,x_0), x_0).$$
\end{Theorem}

\smallskip\noindent Under the assumption of $q$-positivity, it is obvious that D'Angelo's result immediately implies the following:

\medskip
\begin{Corollary}
\label{qposcor}
If $M$ is $q$-positive at $x_0,$ then $$\tilde \Delta_q (M, x_0) = 2 \, {\text gen.val}_{\{w_1, \dots, w_{q-1} \}} \,  \sup_{U \in \, \U(N)} \:\:  \Delta_1 \Big((\I(U,x_0), w_1, \dots, w_{q-1})),x_0\Big).$$
\end{Corollary}

\smallskip\noindent Pseudoconvexity and $q$-positivity relate to each other as follows:

\medskip
\begin{Proposition}[Proposition 2.18  \cite{bazilandreea}]
\label{pscqproppproposition}
If $M$ is pseudoconvex near $x_0$ and $\Delta_q (M, x_0)< + \infty,$ then $M$ and $M_k,$ the hypersurface corresponding to the truncation of order $k$ of the defining function at $x_0,$ are $q$-positive at $x_0$ for all sufficiently large $k.$
\end{Proposition}

Let $V^q$ be the germ of a $q$-dimensional complex variety passing through $x_0.$ Let $G^{n-q+1}$ be the set of all $(n-q+1)$-dimensional complex linear subspaces through $x_0$ in $\C^n.$ Consider the intersection $V^q \cap S$ for $S \in G^{n-q+1}.$ For a generic, namely open and dense, subset $\tilde W$ of $G^{n-q+1},$ $V^q \cap S$ consists of finitely many irreducible one-dimensional components $V^q_{S,k}$ for $k=1, \dots, P.$ We parametrize each such germ of a curve by some open set $U_k \ni 0$ in $\C.$ Thus, $\gamma_S^k : U_k \rightarrow V^q_{S,k},$ where $\gamma_S^k (0)=x_0.$ For every holomorphic germ $f \in \oka_{x_0},$ let $$\tau (f, V^q \cap S) = \max_{k=1, \dots, P} \frac{ {\text ord}_0 \, {\left(\gamma^k_S\right)}^* f}{{\text ord}_0 \, \gamma^k_S}.$$ Likewise, for $r$ the defining function of a real hypersurface $M$ in $\C^n$ passing through $x_0,$ let $$\tau (V^q \cap S, x_0) = \max_{k=1, \dots, P} \frac{ {\text ord}_0 \, {\left(\gamma^k_S\right)}^* r}{{\text ord}_0 \, \gamma^k_S}.$$ In Proposition 3.1 of \cite{catlinsubell}, Catlin showed $\tau (f, V^q \cap S)$ assumes the same value for all $S$ in a generic subset $\tilde W$ of linear subspaces that depends on $V^q,$ so he defined $$\tau(f, V^q) = {\text gen.val}_{S \in \tilde W} \left\{ \tau (f, V^q \cap S) \right\}$$ and $$\tau(\I, V^q) = \min_{f \in \I} \tau(f, V^q).$$ Proposition 3.1 of \cite{catlinsubell} likewise implies that $\tau (V^q \cap S, x_0)$ assumes the same value for all $S$ in the same generic subset $\tilde W$ of linear subspaces depending on $V^q,$ so Catlin defined $$\tau(V^q, x_0) = {\text gen.val}_{S \in \tilde W} \left\{ \tau (V^q \cap S, x_0) \right\}.$$ We will need to explain exactly how the generic subset $\tilde W$ depends on the germ of the $q$-dimensional complex variety $V^q$ when we prove Theorem~\ref{mainthm} in Section~\ref{pfmainthm}, so we defer that discussion.

\medskip
\begin{Definition}
\label{catlinqtypedef}
Let $\I$ be an ideal of holomorphic germs at $x_0,$ then the Catlin $q$-type of the ideal $\I$ is given by $$D_q (\I, x_0)=\sup_{V^q}  \left\{ \tau(\I, V^q) \right\}.$$ 
Let $M$ be a real hypersurface in $\C^n,$ and let $x_0 \in M.$ The Catlin $q$-type of $M$ at $x_0$ is given by $$D_q (M, x_0)=\sup_{V^q}  \left\{ \tau (V^q, x_0) \right\}.$$ In both cases, the supremum is taken over the set of all germs of $q$-dimensional complex varieties $V^q$ passing through $x_0.$ 
\end{Definition}

\noindent Since there is only one $n$-dimensional complex linear subspace passing through $x_0$ in $\C^n,$ $$\Delta_1 (M, x_0) =\tilde \Delta_1 (M, x_0) = D_1 (M , x_0)$$ and $$\Delta_1 (\I, x_0) =\tilde \Delta_1 (\I, x_0) = D_1 (\I , x_0)$$ for any ideal $\I$ in $\oka_{x_0}.$ Therefore, relating these quantities is non-trivial only if $q\geq 2.$

\smallskip We shall also need the following result:

\medskip
\begin{Proposition}[Proposition 4.1 \cite{bazilandreea},\cite{bazilandreeaerr}]
\label{idealleftprop}
Let $\I$ be any ideal in $\oka_{x_0}.$ For any $2 \leq q \leq n,$ $D_q(\I, x_0) \leq \tilde\Delta_q(\I, x_0).$
\end{Proposition}

\section{Proof of results}
\label{pfmainthm}

We claimed in the introduction that $\Delta_q$ and $\tilde \Delta_q$ being simultaneously finite for ideals of germs of holomorphic functions easily follows from D'Angelo's work in \cite{opendangelo}. Here is the proof:

\smallskip\noindent {\bf Proof of Proposition~\ref{normalandtilde}:} Let $\{w_1, \dots, w_{q-1} \}$ be any non-degenerate set of linear forms in $\oka_{x_0}.$ Trivially, 
\begin{equation*}
\begin{split}
\Delta_q(\I, x_0) &= \inf_{\{w_1, \dots, w_{q-1} \}}\:\: \Delta_1\Big( (\I, w_1, \dots, w_{q-1}), x_0 \Big) \\&\leq{\text gen.val}_{\{w_1, \dots, w_{q-1} \}} \:\:\ \Delta_1\Big( (\I, w_1, \dots, w_{q-1}), x_0 \Big) = \tilde \Delta_q(\I, x_0).
\end{split}
\end{equation*}
Applying \eqref{dangeloineq} in Theorem~\ref{opendangelomain} with $q-1$ instead of $q$ yields $${\text gen.val}_{\{w_1, \dots, w_{q-1} \}} \:\:\ \Delta_1\Big( (\I, w_1, \dots, w_{q-1}), x_0 \Big) \leq {\text gen.val}_{\{w_1, \dots, w_{q-1} \}} \:\:\ D \Big( (\I, w_1, \dots, w_{q-1}), x_0 \Big)$$ and $$ \inf_{\{w_1, \dots, w_{q-1} \}} \:\: D \Big( (\I, w_1, \dots, w_{q-1}), x_0 \Big) \leq  \left(\inf_{\{w_1, \dots, w_{q-1} \}} \:\:\Delta_1 \Big( (\I, w_1, \dots, w_{q-1}), x_0 \Big)\right)^{n-q+1},$$ but $$\inf_{\{w_1, \dots, w_{q-1} \}} \:\: D \Big( (\I, w_1, \dots, w_{q-1}), x_0 \Big) = {\text gen.val}_{\{w_1, \dots, w_{q-1} \}} \:\:\ D \Big( (\I, w_1, \dots, w_{q-1}), x_0 \Big)$$ by Proposition~\ref{multgeneric}. \qed

Proving that $\Delta_q$ and $\tilde \Delta_q$ are simultaneously finite when measured at points on the boundary of a domain requires $q$-positivity as well as a slightly more elaborate argument:

\smallskip\noindent {\bf Proof of Proposition~\ref{normalandtildebdry}:} Since the definition of $\Delta_q$ involves taking the infimum over all non-degenerate sets $\{w_1, \dots, w_{q-1} \}$ of linear forms in $\oka_{x_0},$ whereas the definition of $\tilde \Delta_q$ involves taking the generic value over the same sets, it is obvious that $\Delta_q(b \Omega, x_0) \leq \tilde \Delta_q(b \Omega, x_0).$ We thus only have to prove the second inequality. Assume that $\Delta_q(b \Omega, x_0) < +\infty;$ otherwise, $\Delta_q(b \Omega, x_0)=\tilde \Delta_q(b \Omega, x_0)=+\infty,$ and there is nothing to prove. Now, let $k$ be large enough so that $\Delta_q (M, x_0) = \Delta_q (M_k, x_0),$ where $M_k$ is real hypersurface defined by $r_k,$ the polynomial with the same $k$-jet at $x_0$ as the defining function $r$ of $M;$ see equation~\eqref{truncation}. Carry out the polarization of $r_k$ to arrive at the ideal $\I(U,x_0)$ in $\oka_{x_0}$ defined for each unitary matrix $U \in \U(N).$ By Theorem~\ref{dangeloqboundthm}, $$    \sup_{U \in \, \U(N)} \:\: \inf_{\{w_1, \dots, w_{q-1} \}} \:\: D\Big( (\I(U,x_0), w_1, \dots, w_{q-1}), x_0 \Big) \leq  \Big(\Delta_q(M, x_0) \Big)^{n-q}.$$ By Proposition~\ref{multgeneric}, we can substitute the generic value for the infimum as follows:
$$    \sup_{U \in \, \U(N)} \:\: {\text gen.val}_{\{w_1, \dots, w_{q-1}\}} \:\: D\Big( (\I(U,x_0), w_1, \dots, w_{q-1}), x_0 \Big) \leq  \Big(\Delta_q(M, x_0) \Big)^{n-q}$$ Since we are assuming $q$-positivity, Theorem~\ref{opendangelomain} and Corollary~\ref{qposcor} together yield 
\begin{equation*}
\begin{split}
\tilde \Delta_q (M, x_0) &= 2 \, {\text gen.val}_{\{w_1, \dots, w_{q-1} \}} \,  \sup_{U \in \, \U(N)} \:\:  \Delta_1 \Big((\I(U,x_0), w_1, \dots, w_{q-1})),x_0\Big)\\
&\leq 2 \, {\text gen.val}_{\{w_1, \dots, w_{q-1} \}} \,  \sup_{U \in \, \U(N)} \:\: D\Big( (\I(U,x_0), w_1, \dots, w_{q-1}), x_0 \Big).
\end{split}
\end{equation*}
Since $\U(N)$ is compact, Propositions~\ref{multgeneric} and ~\ref{duscprop} imply that the supremum and the generic value can be exchanged, so we obtain that $$ \tilde \Delta_q (M, x_0) \leq 2  \Big(\Delta_q(M, x_0) \Big)^{n-q}$$ by combining the two previous inequalities. If M is pseudoconvex at $x_0,$ by Proposition~\ref{pscqproppproposition}, M is $q$-positive at $x_0,$ so the same inequality must hold. \qed

\medskip
Our proof of the Main Theorem~\ref{mainthm} is motivated by a remark Catlin made on pp.147-8 of \cite{catlinsubell}. Catlin's remark can be paraphrased as follows: In order to compare the D'Angelo $q$-type with his own notion of $q$-type, it would be necessary to piece together one-dimensional varieties in order to get a $q$-dimensional one.  We start with the following lemma, which constructs a very simple $q$-dimensional variety containing a given curve:

\medskip
\begin{Lemma}
\label{cylinderlemma}
Let $2 \leq q \leq n-1.$ Given a holomorphic curve $\Gamma$ passing through a point $x_0 \in \C^n$ and a $(q-1)$-dimensional hyperplane $Z$ passing through $x_0$ and satisfying that the tangent line to the curve $\Gamma$ at $x_0$ is not contained in $Z,$ there exists a germ of a $q$-dimensional cylinder $C^q$ at $x_0$ that contains the curve $\Gamma$ and whose tangent space at $x_0$ contains $Z.$
\end{Lemma}

\begin{Remark}
If $x_0$ is a singular point of the curve $\Gamma,$ the curve could have multiple tangent lines, one for each branch of the curve. In that case, our assumption that the tangent line to the curve $\Gamma$ at $x_0$ is not contained in $Z$ means none of the tangent lines of $\Gamma$ at $x_0$ are contained in $Z.$
\end{Remark}

\smallskip\noindent {\bf Proof:} Let $\Gamma\subset {\mathbb C}^n$ be a holomorphic curve  ($\dim \Gamma =1$) 
given by the equations
$$f_1(z_1,z_2,...,z_n)=0$$

$$f_2(z_1,z_2,...,z_n)=0 $$

$$\vdots$$

$$f_{n-1}(z_1,z_2,...,z_n)=0$$

subject to the condition
$$rank(\partial f_i/\partial z_j)_{1 \leq i \leq n-1;\; 1\leq j \leq n} =n-1$$ that holds generically, namely everywhere except at some isolated points of $\Gamma.$

Let $x_0=(x_1^0, \dots, x_n^0).$ The $(q-1)$-hyperplane $Z$ is given parametrically as

$$x_i=\Sigma_{j=1}^{q-1}u_{ji}t_j + x_i^0,\:\: i=1,2,...,n$$
with the condition
$$rank(u_{ji})=q-1.$$

We have assumed that the tangent line to the curve $\Gamma$ at $x_0$ is not contained in $Z,$ so we can define the cylinder $C^q$ determined by the curve $\Gamma$ and the hyperplane 
$Z$ as follows:

$$C^q:=\{M\in {\mathbb C}^n \; \text{such that there exists}\; P\in \Gamma \;\text{with line}\; MP\; \text{parallel to} \; Z\}.$$ 

By taking the coordinates for $M=(z_1,...,z_n)$ and $P=(x_1,...,x_n),$ we obtain

$$z_i-x_i=\Sigma_{j=1}^{q-1}u_{ji}t_j .$$

From the condition $P\in \Gamma,$ we get the parametric equations of the cylinder $C^q$:

$$f_l(-\Sigma_{j=1}^{q-1}u_{j1}t_j+z_1, ..., -\Sigma_{j=1}^{q-1}u_{jn}t_j+z_n)=0, \:\: l=1,...,n-1.$$

Differentiation with respect to parameters $t_1, \dots, t_{q-1}$ yields

$$\frac{\partial f_l}{\partial t_j}=\sum_{i=1}^n \frac{\partial f_l}{\partial x_i}u_{ji}.$$

The condition that the tangent to the curve $\Gamma$ at $x_0$ is not parallel to $Z$ implies

$$rank \left(\frac{\partial f_l}{\partial t_j}\right)_{1 \leq l \leq n-1;\; 1\leq j \leq q-1}=q-1$$ generically in the neighborhood of $x_0.$

Of course, we have that $\Gamma \subset C^q$ and that the tangent space of $C^q$ at $x_0$ contains $Z$ by construction. \qed

\begin{Remark}
For each of $l=1,...,n-1,$ the hypersurface with equation $f_l(z_1,...,z_n)=0$ contains the curve $\Gamma$.
\end{Remark}

\medskip
\begin{Proposition}
\label{idealprop}
Let $\I$ be an ideal of germs of holomorphic functions at $x_0,$ then for $2 \leq q \leq n-1,$ $$\tilde \Delta_q(\I, x_0) \leq D_q(\I, x_0).$$
\end{Proposition}

\smallskip\noindent {\bf Proof:} We proceed as in the proof of Proposition~\ref{idealleftprop}. Let $D_q(\I, x_0)=t < +\infty,$ else the estimate is trivially true. Assume $\tilde \Delta_q(\I, x_0) >t.$ By definition, $$\tilde \Delta_q(\I, x_0)={\text gen.val}_{\{w_1, \dots, w_{q-1} \}} \:\:\ \Delta_1\Big( (\I, w_1, \dots, w_{q-1}), x_0 \Big),$$ where the generic value is taken over all non-degenerate sets $\{w_1, \dots, w_{q-1} \}$ of linear forms in $\oka_{x_0}.$ Thus, there exist $t'>t$ and a curve $\Gamma$ passing through $x_0$ such that $$\:\: \inf_{g \in (\I, w'_1, \dots, w'_{q-1})} \:\:\frac{ {\text ord}_0 \, \Gamma^* g}{{\text ord}_0 \, \Gamma}=t',$$ for some non-degenerate set $\{w'_1, \dots, w'_{q-1} \}$ of linear forms in $\oka_{x_0}.$ The variety corresponding to the non-degenerate set $\{w'_1, \dots, w'_{q-1} \}$ of linear forms is a $(n-q+1)$-dimensional hyperplane, which we will call $H.$ Clearly, $\Gamma \subset H$ as each $w'_j$ has vanishing order $1$ so any value starting at $2,$ the lower bound for $\tilde \Delta_q(b \Omega, x_0),$ can only be achieved by a curve that sits in the zero set of $\{w'_1, \dots, w'_{q-1} \}.$ Let $Z$ be the $(q-1)$-dimensional hyperplane passing through $x_0$ that is transversal to $H,$ i.e., $\dim_\C (H \oplus Z)=n.$ Since $\Gamma \subset H,$ the tangent lines to the curve $\Gamma$ at $x_0$ (there could be several if the point is singular) are not contained in $Z.$  By Lemma~\ref{cylinderlemma}, there exists a germ of a $q$-dimensional cylinder $C^q_Z$ at $x_0$ that contains the curve $\Gamma$ and whose tangent space at $x_0$ contains $Z.$ We would like to show that $\Gamma$ is one of the intersection curves along which $\tau((\I, w'_1, \dots, w'_{q-1}), C^q_Z)$ is measured by examining Catlin's construction in Proposition 3.1 of \cite{catlinsubell}. 

In the proof of Proposition 3.1 from \cite{catlinsubell}, Catlin removed three different sets $W_1,$ $W_2,$ and $W_3$ from $G^{n-q+1}$ in order to arrive at his generic set $\tilde W$ on which his $q$-type is computed. We wish to show that $H \in \tilde W,$ namely that $H$ cannot belong to any of the three sets $W_1,$ $W_2,$ and $W_3$ taken out from $G^{n-q+1}.$ First, to the germ of the variety $C^q_Z,$ there corresponds the ideal $I=\I(C^q_Z)$ in the ring $\oka_{x_0}$ of all germs of holomorphic functions vanishing on $C^q_Z.$ By a translation, one can suppose $x_0=0$ and thus denote by $\oka$ the ring $\oka_{x_0}.$ Let $_k\oka$ denote the subring of $\oka$ consisting of germs of holomorphic functions of only the first $k$ variables.
%

Now, let $V=C^q_Z,$  and consider as in \cite{catlinsubell} the tangent cone of 
$V$ at the origin, $V'=\{z\in {\mathbb C}^n ; \; \tilde f (z)=0, \; f\in I\}$, where for any function $f$, 
$\tilde f$ is the homogeneous polynomial given by the leading terms in its Taylor expansion at the origin. By \cite{whitney} the tangent cone $V'$ can also be defined by 
\begin{equation}
\label{2nddescr}
V'=\left\{\lim_{|z|\rightarrow 0}\frac{cz}{|z|};\;z\in V, \;c\in {\mathbb C}\right\}
\end{equation}
and $V'$ has dimension $q$ since $V$ has dimension $q$. The curve $\Gamma$ in the $(n-q+1)$-dimensional hyperplane $H$ has as its tangent cone at the origin a union of finitely many lines $L_1, \dots , L_s$ through the origin that are all contained in $H.$ 

We take advantage of the fact that the $q$-dimensional variety $C^q_Z$ has a particularly simple structure by construction.  
Since the hyperplane $H$ is given by a non-degenerate set $w'_1, \dots , w'_{q-1}$ of linear forms in $\oka$, by making a linear change of variables, 
we can assume that the equations of the $(n-q+1)$-dimensional hyperplane $H$ are $z_1=0, \dots , z_{q-1}=0$. Then the equations of the curve $\Gamma$ have the following form (see Lemma~\ref{cylinderlemma}):

$$f_1=z_1=0,\dots , f_{q-1}=z_{q-1}=0, f_q=f_q(z_1,\dots ,z_{q-1}, z_q,\dots , z_n)=0,$$ $$\dots , f_{n-1}=f_{n-1}(z_1,\dots ,z_{q-1}, z_q, \dots , z_n)=0.$$

By choosing the parametric equations of the hyperplane $Z$ to be of the form
$$z_1=t_1, \dots ,z_{q-1}=t_{q-1}, z_q=0, \dots, z_n=0,$$
we get for the parametric equations of the cylinder $C^q_Z$ (see the proof of Lemma~\ref{cylinderlemma}):

$$z_1-t_1=0, \dots , z_{q-1}-t_{q-1}=0, f_q(0,\dots ,0, z_q, \dots , z_n)=0,$$
$$\dots , f_{n-1}(0, \dots ,0, z_q, \dots , z_n)=0.$$

Therefore, the equations for the cylinder $C^q_Z$ are given by

$$f_q=f_q(0,\dots ,0, z_q, \dots , z_n)=0,\dots , f_{n-1}=f_{n-1}(0, \dots ,0, z_q, \dots , z_n)=0,$$
which are the equations of the curve $\Gamma$ in the hyperplane $H$.

Let us denote by $I_1$ the ideal generated by $(f_q, \dots, f_{n-1})$ in the ring 
$\oka$ and by $J$ the ideal generated by $(f_q, \dots, f_{n-1})$ in the ring 
${\oka}_{q,\dots ,n}$, the ring of germs of holomorphic functions at $0$ in the hyperplane $H\cong {\mathbb C}^{n-q+1}$. We have the relations

$$I=rad (I_1), {\oka}_{q,\dots ,n}\subset {\oka}, \;J\subset I.$$
Let $_k{\oka}_{q,\dots ,n}$ denote the subring of germs of holomorphic functions of 
only the variables $z_q, \dots , z_k$, where $k=q+1, \dots , n$. 

Using Gunning's Local Parametrization Theorem from p.15 of \cite{gunning}, since the dimension of $\Gamma$ in $H$ is $1$, one can choose functions 
$g_k, \;k=q+1, \dots , n$ and linear coordinates $(z_q, z_{q+1}, \dots , z_n)$ such that: (1) $g_k\in J\cap_{k-1}{\oka}_{q,\dots ,n} [z_k]$, (2) $g_k$ is a Weierstrass polynomial in $z_k$ of degree $m_k$, and (3) $_q{\oka}_{q,\dots ,n}\cap J=0.$

Moreover, at the last step when $g_{q+1}$ is chosen, by a linear change of variables in $z_q$ and $z_{q+1}$, we can suppose that the $(n-q)$-dimensional hyperplane given by the equation $z_q=0$ in the hyperplane $H$ does not contain the lines $L_1, \dots , L_s$. This can be seen as follows: Take the projective space 
${\mathbb P}^{n-q}$ associated to the hyperplane H, and let $\pi :H\rightarrow {\mathbb P}^{n-q}$ be the canonical projection. The lines $L_1, \dots , L_s$ will define $s$ points in the projective space ${\mathbb P}^{n-q}$ and we can choose a hyperplane $Y$ of dimension $(n-q-1),$ which does not contain any of these points. ${\pi}^{-1}(Y),$ the inverse image of $Y$ in $H,$ then does not contain the lines $L_1, \dots , L_s$. Now, the orthogonal line to ${\pi}^{-1}(Y)$ in $H$ is the $z_q$-axis as $_q{\oka}_{q,\dots ,n}\cap J=0$ by condition (3) above.

We remark that $_{q-1}\oka\cap I=0$. Indeed, if $h(z_1, \dots , z_{q-1})\in I$, then we have for some $j$ that $h^j \in I_1,$ so 
$h^j=a_qf_q+ \dots + a_{n-1}f_{n-1}$, where $a_q, \dots, a_{n-1}$ are germs of holomorphic functions at $0$. Setting in this relation $z_q=0, \dots, z_n=0$ yields
$h^j=0$, hence $h=0$. Since $ {\oka}_{q,\dots ,n}\subset {\oka}, \;J\subset I$, we obtain that the functions $g_k, \; k=q+1, \dots ,n$ satisfy the conditions:  
(1) $g_k\in I\cap_{k-1}\oka [z_k]$, (2) $g_k$ is a Weierstrass polynomial in $z_k$ of degree $m_k$. Furthermore, from the conditions $_{q-1}\oka\cap I=0$ and $_q{\oka}_{q,\dots ,n}\cap J=0,$ it follows that (3) $_q\oka\cap I=0$.

From the existence of $g_k, \; k=q+1, \dots ,n,$ by using some facts about integral extensions (see \cite{gunning}), one can conclude that there exist Weierstrass polynomials $p_k\in I\cap_q\oka [z_k], k=q+1,\dots, n$. By choosing each function $p_k$ to be one of the minimal degree with this property, then since $I$ equals its radical, one may suppose that each $p_k$ is a product of $p_{k,j}, j=1,\dots, N_k$ of distinct irreducible in $_q\oka [z_k]$. It follows that the discriminant $d_k$ of such a product is in the ring $_q\oka$ and not equal to the zero-element (see \cite{catlinsubell}). 

Since the $q$-dimensional variety $V=C^q_Z$ has a particularly simple structure by construction, it follows from the second description of $V'$ (expression \eqref{2nddescr}) that the tangent cone of $V$ is given by the finite union of $q$-dimensional hyperplanes $Z\times L_1, \dots , Z\times L_s$. 

Catlin defined in \cite{catlinsubell} another conic variety $X'\subset V'$ by 

$$X'=\left\{z\in V'; \; z_q{\tilde g}_q\prod_{q+1}^n {\tilde d}_k(z)=0\right\},$$
where  ${\tilde g}_q=1$ in our case since the cylinder variety $V=C^q_Z$ is pure $q$-dimensional. It follows that $X'$ is a conic subvariety of dimension $q-1$ in ${\mathbb C}^n,$ which gives rise to a projective subvariety $\tilde X$ in ${\mathbb P}^{n-1}$ of dimension $q-2$. If $S \in G^{n-q+1}$ and $\tilde S$ denotes the corresponding projective hyperplane of dimension $n-q$ in ${\mathbb P}^{n-1}$, then
the set $W_1$ consists of all $(n-q+1)$-dimensional complex hyperplanes  $S$ such that $\tilde S\cap\tilde X\neq\emptyset$. We have to show that $\tilde H\cap\tilde X=\emptyset$. Firstly, we have $\displaystyle V'\cap H=\cup_{j=1}^s L_j,$ hence $\tilde V\cap\tilde H$ has $s$ points corresponding to the lines $L_1, \dots , L_s$. For the intersection $X'\cap H,$ we have to inspect the equation of $X'$ when $z_1=0, \dots , z_{q-1}=0$. We obtain $z_q\prod_{q+1}^n {\tilde d}_k(0,0, \dots ,0, z_q)=0$ whose only solution is $z_q=0$ because ${\tilde d}_k$ is a homogeneous polynomial for every $k.$ The $(n-q)$-dimensional hyperplane of $H$ given by this solution does not contain any of the lines $L_1, \dots , L_s$ as it is ${\pi}^{-1}(Y)$ described above. It follows that $\tilde H\cap\tilde X=\emptyset$, i.e. $H\notin W_1$.

Secondly, to ensure the intersection $C^q_Z \cap S$ for $S \in G^{n-q+1}$ behaves well, a good notion of transversality has to apply. 
For every $S \in G^{n-q+1},$ there corresponds a projective plane $\tilde S$ of dimension $n-q$ in $\pj^{n-1}.$ Generically, $\tilde C^q_Z \cap \tilde S$ consists of finitely many points $\tilde z^1, \dots, \tilde z^D$ with transverse intersections, meaning that each $\tilde z^i$ is a smooth point of $ \tilde C^q_Z$ such that the tangent spaces satisfy $T_{\tilde z^i} \tilde C^q_Z \cap T_{\tilde z^i} \tilde S = 0$ for $i=1, \dots, D.$ In \cite{catlinsubell} $W_2$ is the subset of $G^{n-q+1}$ where this generic behavior does not take place. As we have seen $\tilde V\cap\tilde H$ has $s$ points corresponding to the lines $L_1, \dots , L_s$. It follows that $D=s,$ and we denote by 
$\tilde z^i$ the point corresponding to $L_i$. Since $L_i$ is a line in $Z\times L_i$ and in $H,$ we have that $T_{\tilde z^i} \tilde V=Z\times L_i/L_i\cong Z$ and $T_{\tilde z^i} \tilde H=H/L_i$. $Z\cap H=0,$ however. It follows that $T_{\tilde z^i} \tilde V \cap T_{\tilde z^i} \tilde H = 0$ for all $i=1,...,s$, i.e. $H\notin W_2$.

Finally, suppose $L\in G^{n-q+1}$ is defined by $\displaystyle L=\left\{z; \sum_{j=1}^na^i_jz_j=0, \; i=1, \dots , q-1\right\}$. Let $\displaystyle W_3=\left\{L; \;\det[a^i_j]_{1\leq i, \;j\leq q-1}=0\right\}$. Since the equations of $H$ are $z_1=0, ..., z_{q-1}=0,$ we have $\det[a^i_j]_{1\leq i, \;j\leq q-1}=1,$ and 
we conclude $H \not\in W_3.$

We have shown that $H \in \tilde W,$ so the curve $\Gamma$ enters into the computation of $\tau((\I, w'_1, \dots, w'_{q-1}), C^q_Z).$ Therefore, $$D_q(\I, x_0) \geq \inf_{g \in (\I, w'_1, \dots, w'_{q-1})} \:\:\frac{ {\text ord}_0 \, \Gamma^* g}{{\text ord}_0 \, \Gamma}=t' > t = D_q(\I, x_0),$$ giving us the needed contradiction. 

\qed

\medskip
\begin{Proposition}
\label{hypersurfaceprop}
Let $\Omega$ in $\C^n$ be a domain with $\smooth$ boundary. Let $x_0 \in b \Omega$ be a point on the boundary of the domain, and let $2 \leq q \leq n-1.$ Then $$ \tilde \Delta_q(b \Omega, x_0)\leq D_q(b \Omega, x_0).$$
\end{Proposition}

\smallskip\noindent {\bf Proof:} The proof is very similar to that of the previous result, Propostion~\ref{idealprop}. Let $D_q(b \Omega, x_0)=t < +\infty,$ else the estimate is trivially true. Assume $\tilde \Delta_q(b \Omega, x_0) >t.$ By definition, $$\tilde \Delta_q(b \Omega, x_0)={\text gen.val}_{\{w_1, \dots, w_{q-1} \}} \:\:\ \Delta_1\Big( (\I(b \Omega), w_1, \dots, w_{q-1}), x_0 \Big),$$ where the generic value is taken over all non-degenerate sets $\{w_1, \dots, w_{q-1} \}$ of linear forms in $\oka_{x_0}$ and $$\Delta_1\Big( (\I(b \Omega), w_1, \dots, w_{q-1}), x_0 \Big) = \sup_{\varphi \in \cv(n,x_0)} \:\: \inf_{g \in(\I(b \Omega), w_1, \dots, w_{q-1})} \:\:\frac{ {\text ord}_0 \, \varphi^* g}{{\text ord}_0 \, \varphi}.$$ Therefore, there exist $t'>t$ and a curve $\Gamma$ passing through $x_0$ such that $$\:\: \inf_{g \in (\I(b \Omega), w_1, \dots, w_{q-1})} \:\:\frac{ {\text ord}_0 \, \Gamma^* g}{{\text ord}_0 \, \Gamma}=t',$$ for some non-degenerate set $\{w'_1, \dots, w'_{q-1} \}$ of linear forms in $\oka_{x_0}.$ We will call $H$ the variety corresponding to the non-degenerate set $\{w'_1, \dots, w'_{q-1} \}$ of linear forms, which is a $(n-q+1)$-dimensional hyperplane. Clearly, $\Gamma \subset H$ as each $w'_j$ has vanishing order $1$ so any value starting at $2,$ the lower bound for $\tilde \Delta_q(b \Omega, x_0),$ can only be achieved by a curve that sits in the zero set of $\{w'_1, \dots, w'_{q-1} \}.$ Therefore, $\inf_{g \in (\I(b \Omega), w_1, \dots, w_{q-1})}$ is achieved for $g=r,$ the defining function of the domain $\Omega.$ Let $Z$ be the $(q-1)$-dimensional hyperplane passing through $x_0$ that is transversal to $H,$ i.e., $\dim_\C (H \oplus Z)=n.$ Since $\Gamma \subset H,$ the tangent line to the curve $\Gamma$ at $x_0$ is not contained in $Z.$  By Lemma~\ref{cylinderlemma}, there exists a germ of a $q$-dimensional cylinder $C^q_Z$ at $x_0$ that contains the curve $\Gamma$ and whose tangent space at $x_0$ contains $Z.$ By the same analysis as in the proof of Propostion~\ref{idealprop}, $H \in \tilde W,$ so the curve $\Gamma$ enters into the computation of $\tau(C^q_Z, x_0).$ Thus, $$D_q(b \Omega, x_0) \geq \frac{ {\text ord}_0 \, \Gamma^* r}{{\text ord}_0 \, \Gamma}=t' > t = D_q(\I, x_0),$$ and we have obtained the contradiction we sought. 
  \qed

\medskip\noindent {\bf Proof of Theorem~\ref{mainthm}:} First, we prove part (i). For q=n Proposition~\ref{normalandtilde} along with Theorem ~\ref{bazilandreeathm1.1} yield that $\Delta_n(\I,x_0)=\tilde \Delta_n(\I,x_0)=D_n(\I,x_0).$ As a result, the equality only needs to be proven for $2 \leq q \leq n-1.$ The inequality $D_q(\I, x_0) \leq \tilde \Delta_q(\I, x_0)$ is a consequence of Theorem~\ref{bazilandreeathm1.1}, while $\tilde \Delta_q(\I, x_0) \leq D_q(\I, x_0)$ follows from Proposition~\ref{idealprop}.

To prove part (ii), we note that $D_q(b \Omega, x_0) \leq \tilde \Delta_q(b \Omega, x_0)$ is a consequence of Theorem~\ref{bazilandreeathm1.2}, while the reverse inequality follows from Proposition~\ref{hypersurfaceprop}.  \qed

\medskip\noindent {\bf Proof of Corollary~\ref{maincor}:} (i) We combine Proposition~\ref{normalandtilde} with Theorem~\ref{mainthm} part (i).

\smallskip\noindent (ii) Proposition~\ref{normalandtildebdry} and Theorem~\ref{mainthm} part (ii) together give the result. \qed

\bibliographystyle{plain}
\bibliography{CatlinDAngeloType}

\begin{thebibliography}{10}

\bibitem{bazilandreea}
Vasile Brinzanescu and Andreea~C. Nicoara.
\newblock On the relationship between {D}'{A}ngelo {$q$}-type and {C}atlin
  {$q$}-type.
\newblock {\em J. Geom. Anal.}, 25(3):1701--1719, 2015.

\bibitem{bazilandreeaerr}
Vasile Brinzanescu and Andreea~C. Nicoara.
\newblock Correction to: {O}n the relationship between {D}'{A}ngelo {$q$}-type
  and {C}atlin {$q$}-type.
\newblock {\em J. Geom. Anal.}, 30(1):1171--1172, 2020.

\bibitem{catlinsubell}
David Catlin.
\newblock Subelliptic estimates for the {$\overline\partial$}-{N}eumann problem
  on pseudoconvex domains.
\newblock {\em Ann. of Math. (2)}, 126(1):131--191, 1987.

\bibitem{opendangelo}
John~P. D'Angelo.
\newblock Real hypersurfaces, orders of contact, and applications.
\newblock {\em Ann. of Math. (2)}, 115(3):615--637, 1982.

\bibitem{dangeloftc}
John~P. D'Angelo.
\newblock Finite-type conditions for real hypersurfaces in {${\bf C}^n$}.
\newblock In {\em Complex analysis ({U}niversity {P}ark, {P}a., 1986)}, volume
  1268 of {\em Lecture Notes in Math.}, pages 83--102. Springer, Berlin, 1987.

\bibitem{fassina}
Martino Fassina.
\newblock A {R}emark on {T}wo {N}otions of {O}rder of {C}ontact.
\newblock {\em J. Geom. Anal.}, 29(1):707--716, 2019.

\bibitem{gunning}
R.~C. Gunning.
\newblock {\em Lectures on complex analytic varieties: {T}he local
  parametrization theorem}.
\newblock Mathematical Notes. Princeton University Press, Princeton, N.J.,
  1970.

\bibitem{kohnacta}
J.~J. Kohn.
\newblock Subellipticity of the {$\bar \partial $}-{N}eumann problem on
  pseudo-convex domains: sufficient conditions.
\newblock {\em Acta Math.}, 142(1-2):79--122, 1979.

\bibitem{tougeron}
Jean-Claude Tougeron.
\newblock {\em Id\'eaux de fonctions diff\'erentiables}.
\newblock Springer-Verlag, Berlin, 1972.
\newblock Ergebnisse der Mathematik und ihrer Grenzgebiete, Band 71.

\bibitem{whitney}
Hassler Whitney.
\newblock {\em Complex analytic varieties}.
\newblock Addison-Wesley Publishing Co., Reading, Mass.-London-Don Mills, Ont.,
  1972.

\end{thebibliography}

\end{document}